\documentstyle[11pt]{article}
\newcount\VOL\VOL=1\newcount\YEAR\YEAR=1996
\def\firstpage{1}\def\lastpage{1000}
\def\received{}\def\revised{}
\def\communicated{}
\font\eightrm=cmr8
\font\caps=cmcsc10                    % Theorem, Lemma etc
\font\Caps=cccsc10 scaled \magstep1   % Title

\makeatletter
\@specialpagefalse\headheight=8.5pt
\def\DocMath{}
\renewcommand{\@evenhead}{%
    \ifnum\thepage>\lastpage\rlap{\thepage}\hfill%
    \else\rlap{\thepage}\slshape\leftmark\hfill\caps\SAuthor\hfill\fi}%
\renewcommand{\@oddhead}{%
    \ifnum\thepage=\firstpage{\DocMath\hfill\llap{\thepage}}%
    \else{\slshape\rightmark}\hfill\caps\STitle\hfill\llap{\thepage}\fi}%
\makeatother

\def\TSkip{\bigskip}
\newbox\TheTitle{\obeylines\gdef\GetTitle #1
\ShortTitle  #2
\SubTitle    #3
\Author      #4
\ShortAuthor #5
\EndTitle
{\setbox\TheTitle=\vbox{\baselineskip=20pt\let\par=\cr\obeylines%
\halign{\centerline{\Caps##}\cr\noalign{\medskip}\cr#1\cr}}%
	\copy\TheTitle\TSkip\TSkip%
\def\next{#2}\ifx\next\empty\gdef\STitle{#1}\else\gdef\STitle{#2}\fi%
\def\next{#3}\ifx\next\empty%
    \else\setbox\TheTitle=\vbox{\baselineskip=20pt\let\par=\cr\obeylines%
    \halign{\centerline{\caps##} #3\cr}}\copy\TheTitle\TSkip\TSkip\fi%
\centerline{\caps #4}\TSkip\TSkip%
\def\next{#5}\ifx\next\empty\gdef\SAuthor{#4}\else\gdef\SAuthor{#5}\fi%
\ifx\received\empty\relax
    \else\centerline{\eightrm Received: \received}\fi%
\ifx\revised\empty\TSkip%
    \else\centerline{\eightrm Revised: \revised}\TSkip\fi%
\ifx\communicated\empty\relax
    \else\centerline{\eightrm Communicated by \communicated}\fi\TSkip\TSkip%
\catcode'015=5}}\def\Title{\obeylines\GetTitle}
\def\Abstract{\begingroup\narrower
    \parskip=\medskipamount\parindent=0pt{\caps Abstract. }}
\def\EndAbstract{\par\endgroup\TSkip\TSkip}
\newbox\TheAdd\def\Addresses{\vfill\copy\TheAdd\vfill
    \ifodd\number\lastpage\vfill\eject\phantom{.}\vfill\eject\fi}
{\obeylines\gdef\GetAddress #1
\Address #2 
\Address #3
\Address #4
\EndAddress
{\def\xs{6truecm}
\setbox0=\vtop{{\obeylines\hsize=\xs#1}}\def\next{#2}
\ifx\next\empty % 1 address
     \setbox\TheAdd=\hbox to\hsize{\hfill\copy0\hfill}
\else\setbox1=\vtop{{\obeylines\hsize=\xs#2}}\def\next{#3}
\ifx\next\empty % 2 addresses
     \setbox\TheAdd=\hbox to\hsize{\hfill\copy0\hfill\copy1\hfill}
\else\setbox2=\vtop{{\obeylines\hsize=\xs#3}}\def\next{#4}
\ifx\next\empty\ % 3 addresses
     \setbox\TheAdd=\vtop{\hbox to\hsize{\hfill\copy0\hfill\copy1\hfill}
	        \vskip20pt\hbox to\hsize{\hfill\copy2\hfill}}
\else\setbox3=\vtop{{\obeylines\hsize=\xs#4}}
     \setbox\TheAdd=\vtop{\hbox to\hsize{\hfill\copy0\hfill\copy1\hfill}
	        \vskip20pt\hbox to\hsize{\hfill\copy2\hfill\copy3\hfill}}
\fi\fi\fi\catcode'015=5}}\gdef\Address{\obeylines\GetAddress}

\begin{document} 
%%%% ------------- fill in your data below this line  -------------------

\Title{Geometric zeta-functions of locally symmetric spaces}
\ShortTitle{Geometric zeta-functions}
\SubTitle   
\Author{Anton Deitmar} 
\ShortAuthor{A. Deitmar}
\EndTitle

\Abstract{The theory of geometric zeta functions for locally symmetric spaces as initialized by Selberg and continued by numerous mathematicians is generalized to the case of higher rank spaces. We show analytic continuation, describe the divisor in terms of tangential cohomology and in terms of group cohomology which generalizes the Patterson conjecture. We also extend the range of zeta functions in considering higher dimensional flats.}

1991 Mathematics Subject Classification: Primary 11F72, Secondary: 11M41, 22E40, 43A85, 53C35, 58F20.
\EndAbstract

\Address{Anton Deitmar\\ Math. Inst. d. Universit\"at\\ INF 288\\ 69126 Heidelberg\\ Germany}
\Address{}
\Address{}
\Address{}
\EndAddress
%
% Make sure the last tex command in your manuscript
% before the first \end or \bye is the command  \Addresses
%
%---------------------Here the prologue ends---------------------------------
%--------------------Here the manuscript starts------------------------------

\input amssym.def
\input amssym

\tableofcontents

\def \a{{{\frak a}}}
\def \al{\alpha}
\def \ar{{\alpha_r}}
\def \A{{\Bbb A}}
\def \Ad{{\rm Ad}}
\def \b{{{\frak b}}}
\def \bs{\backslash}
\def \B{{\cal B}}
\def \cent{{\rm cent}}
\def \C{{\bf C}}
\def \CA{{\cal A}}
\def \CB{{\cal B}}
\def \CE{{\cal E}}
\def \CF{{\cal F}}
\def \CG{{\cal G}}
\def \CH{{\cal H}}
\def \CM{{\cal M}}
\def \CN{{\cal N}}
\def \CP{{\cal P}}
\def \CQ{{\cal Q}}
\def \CO{{\cal O}}
\def \det{{\rm det}}
\def \End{{\rm End}}
\def \Fx{{\frak x}}
\def \FX{{\frak X}}
\def \g{{{\frak g}}}
\def \ga{\gamma}
\def \Ga{\Gamma}
\def \h{{{\frak h}}}
\def \Hom{{\rm Hom}}
\def \Im{{\rm Im}}
\def \Ind{{\rm Ind}}
\def \k{{{\frak k}}}
\def \K{{\cal K}}
\def \la{\lambda}
\def \lap{\triangle}
\def \La{\Lambda}
\def \m{{{\frak m}}}
\def \mod{{\rm mod}}
\def \n{{{\frak n}}}
\def \name{\bf}
\def \N{\Bbb N}
\def \ord{{\rm ord}}
\def \O{{\cal O}}
\def \p{{{\frak p}}}
\def \ph{\varphi}
\def \prf{{\bf Proof: }}
\def \qed{\hfill {$\Box$} 

$ $

}
\def \Q{\Bbb Q}
\def \res{{\rm res}}
\def \R{{\Bbb R}}
\def \Re{{\rm Re \hspace{1pt}}}
\def \ra{\rightarrow}
\def \rank{{\rm rank}}
\def \supp{{\rm supp}}
\def \t{{{\frak t}}}
\def \T{{\Bbb T}}
\def \tr{{\hspace{1pt}\rm tr\hspace{2pt}}}
\def \vol{{\rm vol}}
\def \V{{\cal V}}
\def \z{\zeta}
\def \Z{\Bbb Z}
\def \={\ =\ }

\newcommand{\rez}[1]{\frac{1}{#1}}
\newcommand{\der}[1]{\frac{\partial}{\partial #1}}
\newcommand{\binom}[2]{\left( \begin{array}{c}#1\\#2\end{array}\right)}

\newcounter{lemma}
\newcounter{corollary}
\newcounter{proposition}
\newcounter{theorem}

\newtheorem{conjecture}{\stepcounter{lemma} \stepcounter{corollary} 	
	\stepcounter{proposition}\stepcounter{theorem}Conjecture}[section]
\newtheorem{lemma}{\stepcounter{conjecture}\stepcounter{corollary}	
	\stepcounter{proposition}\stepcounter{theorem}Lemma}[section]
\newtheorem{corollary}{\stepcounter{conjecture}\stepcounter{lemma}
	\stepcounter{proposition}\stepcounter{theorem}Corollary}[section]
\newtheorem{proposition}{\stepcounter{conjecture}\stepcounter{lemma}
	\stepcounter{corollary}\stepcounter{theorem}Proposition}[section]
\newtheorem{theorem}{\stepcounter{conjecture} \stepcounter{lemma}
	\stepcounter{corollary}	\stepcounter{proposition}Theorem}[section]

$ $

\begin{center} {\bf Introduction} \end{center}

In his celebrated paper on the trace formula 
\cite{Sel} A. Selberg introduced the zeta function
$$
Z(s) = \prod_c \prod_{N\geq 0} (1-e^{-(s+N)l(c)}),
$$
where the first product is taken over all primitive closed geodesics in a compact Riemannian surface of genus $\geq 2$ and $l(c)$ denotes the length of the geodesic $c$.
Selberg proved that $Z(s)$ extends to an entire function on $\C$, that it satisfies a functional equation as $s$ is replaced by $1-s$ and that 
a generalized Riemann hypothesis holds for the function $Z$.
As the trace formula parallels the Poisson summation formula the zeta function parallels the Riemann zeta function up to the difference that one can prove the generalized Riemann hypothesis for the latter.
Many authors have studied this zeta function more thoroughly, see \cite{BuOl-buch}, \cite{CaVo}, 
    \cite{Efrat-det}, \cite{Efrat-dyn}, \cite{Hej}, \cite{Sarnak}, 
    \cite{Venkov}, \cite{Voros} or \cite{Els} and the references therein.

The trace formula has been generalized to semisimple Lie groups \cite{Wall-stf} and to reductive groups over $\Q$ \cite{Art-1}, \cite{Art-2}.
The generalization of the zeta function took somewhat more time. R. Gangolli \cite{Gang} defined a zeta function of Selberg type for compact locally symmetric spaces whose universal covering is a noncompact locally symmetric space of rank one.
M. Wakayama \cite{Wak} extended the theory of Gangolli to allow twists by homogeneous vector bundles.
D. Fried \cite{Fr-at} found a simpler way to get analytic continuation by heat kernel methods. This also allowed him to interprete the zeta function as a determinant of the Laplacian. This differential operator approach is worked out for the rank one case in the book \cite{BuOl-buch}.
A. Juhl \cite{Ju} gave a geometric way to achieve the analytic continuation 
which describes the divisor of the zeta function in terms of Lie algebra 
cohomology.

The case of higher rank seemed impenetrable until H. Moscovici and R. Stanton 
\cite{MS-tors} used supersymmetry arguments to compute traces of certain 
linear combinations of heat operators. 
This made it possible for them to 
get the continuation of the Ruelle zeta function in the cases $SL_3(\R)$ and $SO(p,q)$ with $pq$ odd. The Ruelle zeta function is a rational 
function in some Selberg zeta functions. Their method was generalized in 
\cite{D-Hitors} to compact locally symmetric manifolds of odd dimension.

In this paper we present a geometric approach which works for arbitrary 
locally symmetric manifolds.
The central idea is to use the rank one approach of A. Juhl to shift
down to a Levi component of splitrank one less.
Then one plugs in an Euler-Poincar\'e function to single out the elliptic
classes in that particular Levi component.

We express the vanishing order of the zeta function in terms of nilpotent Lie algebra cohomology, thereby proving the generalized Riemann hypothesis for the zeta function and a determinant formula similar to the determinant formula of Deligne expressing the Hasse-Weil zeta function as an alternating product of determinants of the Frobenius-action on \'etale cohomology. 
We further relate the Lie algebra cohomology to the cohomology of the fundamental group thus proving a version of the Patterson conjecture in our setting.
In the rank one case this conjecture was proven by U. Bunke and M. Olbrich \cite{BuOl}. 

The zeta functions of this paper can also be interpreted as dynamical zeta functions.
In this interpretation they represent new objects in several aspects: In fundamental rank zero they give zeta functions attached to non Anosov dynamical systems and in the case of arbitrary fundamental rank the $L$-functions are attached to dynamical systems of higher dimensional time.

I thank Andreas Juhl for his comments and Ulrich Bunke and Martin Olbrich for communicating their work to me.

$ $

{\bf Notation}

We will write $\N , \Z ,\Q ,\R ,\C$ for the natural, integer, rational, real and complex numbers.

For any locally compact group $G$ we write $\hat{G}$ for the unitary dual, i.e. the set of equivalence classes of irreducible unitary representations.

Throughout we will use small german letters to denote Lie algebras of the corresponding Lie groups which will be denoted by capital roman letters. A subscript zero will indicate the Lie algebra over $\R$, otherwise it will be its complexification, So for example $G$ a Lie group, then we write $\g_0=Lie_\R(G)$ and $\g =\g_0\otimes_\R \C$.

A {\bf virtual vector space} $V$ will be the formal difference of two vector spaces, i.e.:
$$
V = V^+ - V^-.
$$
An endomorphism $A$ of a virtual space $V$ is a pair of endomorphisms $ A^\pm $ on $ V^+ $ and $ V^- $ resp. The {\bf trace} and {\bf determinant} of $A$ are then
$$
\rm{tr} (A) = \rm{tr} (A^+) - \rm{tr} (A^-),\ \ \ \ \ \rm{\rm{det}} (A) = \rm{\rm{det}} (A^+) / \rm{\rm{det}} (A^-).
$$

Every vector space $V$ with $ \Z $-grading will naturally be considered as a virtual vector space by
$$
V^+ = V_{even},\ \ \ \ \ \ \ \ \ \ V^- = V_{odd}.
$$

As an example consider the exterior algebra over the finite dimensional space $V$ then for any endomorphism $A$ of $V$ we have the well known formula
$$
\rm{\rm{det}}( 1 - A ) = \rm{tr} (A \mid \wedge^* V ).
$$

For an endomorphism $A$ of a finite dimensional vector space $V$ denote by
$V_0 = \cup_{n\in \N} \ker A^n$ the generalized kernel of $A$. Since $V_0$ is stable under $A$ it induces an endomorphism of $V' = V/V_0$, which is injective.
Define the {\bf essential determinant} of $A$ to be
$$
\det'(A) = \det(A').
$$

Although our results do not depend on Haar-measure normalizations we will need those for the computations along the way.
All Haar measures occurring are normalized as in \cite{HC-HA1}.
We will make this a bit more precise.
Fix a semisimple connected Lie group $G$ without center and such that $G$ admits a compact Cartan-subgroup.
Recall from \cite{HC-HA1} that the Haar-measure normalization depends on the choice of an invariant bilinear form $B$ on the Lie algebra of $G$.
In the following fix a maximal compact subgroup $K$ of $G$ and write $X$ for the quotient space $G/K$.
We will choose $B$ arbitrary as long as $G$ has no compact Cartan subgroup.
When $G$ has a compact Cartan, however, we will choose $B$ such that the Haar measure it induces equals the absolute value of the {\bf Euler-Poincar\'e measure}, which is the unique Haar measure on $G$ such that for all torsion-free cocompact discrete subgroups $\Ga$ we have
$$
\vol(\Ga \bs G) = (-1)^{\dim X/2} \chi(\Ga \bs X),
$$
here $\chi$ denotes the Euler-characteristic.
Having chosen $B$ we use it to define an invariant metric on the manifold $X$.
This makes $X$ a globally symmetric space.

$ $

\section{Euler-Poincar\'e functions}
In this section we prepare some material which we later will apply to the semi-simple part of a Levi-component. So let $G$ denote a semisimple real reductive group of inner type \cite{Wall-rg1} and fix a maximal compact subgroup $K$.
Let $(\tau ,V_\tau)$ be a finite dimensional unitary representation of $K$ and write $(\breve{\tau},V_{\breve{\tau}})$ for the dual representation. 
Assume that $G$ has a compact Cartan subgroup $T \subset K$.
Let $\g_0 = \k_0 \oplus \p_0$ be the polar decomposition of the real Lie algebra $\g_0$ of $G$ and write $\g = \k +\p$ for its complexification.
Let $\t$ be the complexified Lie algebra of the cartan subgroup $T$.
Choose an ordering of the roots $\Phi(\g ,\t)$ of the pair $(\g ,\t)$. This choice induces a decomposition $\p = \p_- \oplus \p_+$.

\begin{theorem} \label{existf}
For $(\tau ,V_\tau)$ a finite dimensional representation of $K$ there is a compactly supported smooth function $f_\tau$ on $G$ such that for every irreducible unitary representation $(\pi ,V_\pi)$ of $G$ it holds:
$$
\tr\ \pi (f_\tau) \= \sum_{p=0}^{\dim (\p)} (-1)^p \dim (V_\pi \otimes \wedge^p\p \otimes V_{\tau})^K.
$$
\end{theorem}

\prf 
Let $S:=\wedge^*\p_-$. Replacing $G$ by a twofold cover if necessary assume that the adjoint homomorphism $K\ra SO(\p)$ factors over the spin group $Spin(\p)$.
Then K acts via the spin representation on $S^+ :=\wedge^{even} \p_-$ and $S^- :=\wedge^{odd}\p_-$ in a way that the action on $S\otimes S$ coincides with the adjoint action.
From \cite{Lab} we take

\begin{lemma}
For any finite dimensional representation $(\tau ,V_\tau)$ of $K$ there is a compactly supported smooth function $g_\tau$ on $G$ such that for every irreducible unitary representation $(\pi ,V_\pi)$ of $G$ we have
$$
\tr\ \pi (g_\tau) \= \dim (V_\pi \otimes S^+ \otimes V_{\tau})^K
- \dim (V_\pi \otimes S^- \otimes V_{\tau})^K.
$$ 
\end{lemma}\qed

We replace $\tau$ in the lemma with the virtual representation on
$$
(S^+-S^-)\otimes V_{{\tau}}
$$
and we get the desired function at least on a twofold cover $\tilde{G}$ of $G$.
Let the kernel of $\tilde{G}\ra G$ be $\{1,\alpha\}$ then, if necessary, replace $f_\tau$ by $(f_\tau(x) +f_\tau(x\alpha))/2$.
The latter function then lives on $G$ and satisfies the theorem.
\qed

We want to show that $\tr\pi(f_\tau)$ vanishes for a principal series representation $\pi$.
To this end let $P=MAN$ be a parabolic subgroup with $A\subset \exp(\p_0)$.
Let $(\xi ,V_\xi)$ denote an irreducible unitary representation of $M$ and $e^\nu$ a quasicharacter of $A$.
Let $\pi_{\xi ,\nu}:= {\rm Ind}_P^G (\xi \otimes e^{\nu}\otimes 1)$.

\begin{lemma} \label{pivonggleichnull}
We have $\tr\pi_{\xi ,\nu}(f_\tau) =0$.
\end{lemma}

\prf
By Frobenius reciprocity we have for any irreducible unitary representation $\ga$ of $K$:
$$
\Hom_K(\ga ,\pi_{\xi ,\nu}|_K) \cong \Hom_{K_M}(\ga |_{K_M},\xi ),
$$
where $K_M := K\cap M$. 
This implies that $\tr\pi_{\xi ,\nu}(f_\tau)$ does not depend on $\nu$.
Write $\la_\pi$ for the infinitesimal character of the representation $\pi$.
By \cite{AtSch}, (4.20) there is a finite set $E$ of infinitesimal characters such that $\la_\pi\notin E$ implies $\tr\pi(f_\tau) =0$.
Since by varying $\nu$ the set of infinitesimal characters of the representations $\pi_{\xi ,\nu}$ is infinite there is one $\nu$ with $\tr\pi_{\xi ,\nu}(f_\tau) =0$. Since it was independent of $\nu$ this assertion is true for all $\nu$.
\qed

We want to compute the orbital integrals
$$
\O_g(f_\tau) \= \int_{G/G_g}f_\tau(xgx^{-1}) dx.
$$

Recall that an element $g$ of $G$ is called {\bf elliptic} if it lies in a compact Cartan subgroup.

\begin{proposition} \label{orbitalint}
Let $g$ be a semisimple element of the group $G$. If $g$ is not elliptic, then the orbital integral $\O_g(f_\tau)$ vanishes. If $g$ is elliptic we may assume $g\in T$, where $T$ is a Cartan in $K$ and then we have
$$
\O_g(f_\tau) \= \frac{{\tr\ \tau(g)}|W(\t ,\g_g)| \prod_{\alpha \in \Phi_g^+}(\rho_g ,\alpha)}{[G_g:G_g^0]c_g},
$$
where $c_g$ is Harish-Chandra's constant, it does only depend on the centralizer $G_g$ of g. Its value is given for example in \cite{D-Hitors}.
\end{proposition}

\prf
The vanishing of $\O_g(f_\tau)$ for nonelliptic semisimple $g$ is immediate by the lemma above and Harish-Chandra's formula for the Fourier transform of orbital integrals \cite{HC-S}. Now let $g\in T\cap G'$, where $G'$ denotes the set of regular elements. 
At first we consider the above function $g_\tau$.
Assume $\tau$ irreducible and for $t>0$ let $h_t$ denote the "heat trace" attached to $\tau$ as in \cite{BM} (3.2) p. 164.
By \cite{BM} (3.5) p. 165 it follows that $\tr \pi(h_t)=\tr\pi(g_\tau)$ for all $\pi\in\hat{G}$.
Hence also the orbital integrals of $h_t$ and $g_\tau$ coincide.
The orbital integrals of $h_t$ are computed in \cite{BM} p.173.
After applying the Weyl character formula this result becomes
$$
\CO_g(h_t)=\frac{\tr\tau(g)}{g^{\rho -\rho_c}\prod_{\beta\in\Phi^+-\Phi^+_c}(1-g^{-\beta})},
$$
where $\Phi^+$, resp. $\Phi^+_c$ denote the sets of positive, resp compact positive roots and $\rho$ resp. $\rho_c$ as usual the half of the sum of the elements of the latter.
The denominator of this expression also makes sense for $g\in Spin(\p)$ as long as $g$ lies in the maximal torus $\tilde{T}\supset image(T)$, whose weight space decomposition gives rise to the decomposition $\p=\p_+\oplus\p_-$.
For $g\in\tilde{T}$ we have
$$
g^{\rho -\rho_c}\prod_{\beta\in\Phi^+-\Phi^+_c}(1-g^{-\beta}) = \tr(g|S^+-S^-).
$$
So we get
$$
\CO_g(g_\tau)=\CO_g(h_t)=\frac{\tr\tau(g)}{\tr(g|S^+-S^-)}
$$
for $g\in G'\cap T$.
Switching to $f_\tau$ means that one replaces $\tau$ by $\tau\otimes(S^+-S^-)$,
hence we infer
$$
\O_g(f_\tau) \= {\tr\ \tau(g)}.
$$
This proves the proposition in the regular case. 
The general case is derived from this by standard considerations (see \cite{HC-DS}, p32 ff. or \cite{BM}p.173).
\qed

\begin{corollary}
If $G$ is a product of real rank one groups then the orbital integral $\O_g(f_\tau)$ vanishes also for $g$ non-semisimple.
\end{corollary}

\prf
It suffices to consider the case $\rank_\R(G)=1$.
Consider a non-semisimple element $g$. By \cite{Barb}, sec. 6 we get a curve $t\mapsto z_t$ of semisimple elements and a natural number $m$ such that $\O_g(f_\tau)=\lim_{t\rightarrow 0} t^{m/2}\O_{z_t}(f_\tau)$.
The function $x\mapsto \O_x(f_\tau)$ is bounded on semisimple elements by the proposition and therefore we have $\O_g(f_\tau) =0$ for $g$ non semisimple.
\qed

\begin{conjecture} \label{all_orbitalint_vanish}
The above corollary holds without the assumption on $G$.
\qed
\end{conjecture}

For a finite dimensional representation $(\sigma ,V_\sigma)$ of $G$ the Theorem \ref{existf} applies to $\sigma |_K$ to give a function $f_\sigma$.

\begin{proposition}
For the function $f_{{\sigma}}$ we have for any $\pi \in \hat{G}$:
$$
\tr \ \pi(f_{{\sigma}}) \= \sum_{p=0}^{\dim \ \g /\k}(-1)^p \dim \ 
{\rm Ext}_{(\g ,K)}^p (V_{\breve{\sigma}} ,V_\pi),
$$
i.e. $f_{{\sigma}}$ gives the Euler-Poincar\'e numbers of the $(\g ,K)$-modules $(V_{\breve{\sigma}} ,V_\pi)$, this justifies the name Euler-Poincar\'e function.
\end{proposition}

\prf
By definition it is clear that
$$
\tr \ \pi (f_{{\sigma}}) \= \sum_{p=0}^{\dim \ \p} (-1)^p \dim \ H^p(\g ,K,V_{{\sigma}} \otimes V_\pi).
$$
The claim now follows from \cite{BorWall}, p. 52. \qed

\section{The construction of the zeta function} \label{constzeta}

We are going to consider locally symmetric spaces of the type $X_\Ga = \Ga \bs G/K$, where $G$ is a connected semisimple Lie-group without center, $K$ a maximal compact subgroup of $G$ and $\Ga$ a torsion-free cocompact lattice.

As the space $X=G/K$ is contractible it follows that 
$\Ga$ is the fundamental group of $X_\Ga$ and thus there is a bijection 
between the set of nontrivial conjugacy classes of $\Ga$ and the set of 
free homotopy classes of closed paths in $X_\Ga$. 
Each such class $[\ga]$ contains closed geodesics, the union of which form 
a submanifold $X_\ga$ which again is a locally symmetric manifold \cite{DKV}.
Let $h^j(X_\ga)$ denote the j-th Betti number of $X_\ga$ then for $r\ge 0$ we define the higher Euler number of $X_\ga$ as
$$
\chi_{_r}(X_\ga) := \sum_{j=0}^{\dim X_\ga} (-1)^{j+r} \binom{j}{r} h^j(X_\ga),
$$
(compare \cite{D-Prod}). 

Let $\theta$ denote the Cartan involution on $G$ given by the choice of $K$. Fix a $\theta$-stable Cartan subgroup $H$ of split rank 1. Note that such a $H$ doesn't always exist.
It exists only if the absolute ranks of $G$ and $K$ satisfy the relation:
$$
{\rm rank}\ G -{\rm rank}\ K \leq 1.
$$
This certainly holds if $G$ has a compact Cartan subgroup or if the real 
rank of $G$ is one. In the case ${\rm rank}\ G ={\rm rank}\ K$, i.e. 
if $G$ has a compact Cartan there will in general be several $G$-conjugacy 
classes of split rank one Cartan subgroups. 
In the case ${\rm rank}\ G -{\rm rank}\ K = 1$, however, there will only 
be one. 
The number $FR(G):={\rm rank}\ G -{\rm rank}\ K $ is called the 
{\bf fundamental rank} of $X$ or $X_\Ga$.

Write $H=AB$ where $A$ is the connected split component and $B\subset K$ is compact. Choose a parabolic subgroup $P$ with Langlands decomposition $P=MAN$ where $A\subset \exp(\p_0)$. 
Then $K_M =K\cap M$ is a maximal compact subgroup of $M$. 
Fix a finite dimensional representation $(\tau ,V_\tau)$ of $K_M$.
Let $G$ act on itself by conjugation, write $g.x = gxg^{-1}$, write $G.x$ for the orbit, so $G.x = \{ gxg^{-1} | g\in G \}$ as well as $G.S = \{ gsg^{-1} | s\in S , g\in G \}$ for any subset $S$ of $G$.
We are going to consider functions that are supported on the closure of the set $G.(MA)$. Now let $f_\tau$ be as in the preceding proposition but with $M$ taking the place of $G$. 
The choice of the parabolic $P$ induces an ordering on the roots $\Phi(\g ,\a)$, where $\a$ is the complexified Lie-algebra of $A$.
Let $\a_0^+$ denote the positive Weyl chamber and $A^+ := \exp (\a_0^+)$.

For $g\in G$ and $V$ any complex vector space on which $g$ acts linearly let $E(g|V)\subset \R^*_+$ be defined by
$$
E(g|V) \ :=\ \{ |\mu | : \mu {\rm \ is\ an\ eigenvalue\ of}\ g {\rm \ on\ } V\}.
$$
Let $\la_{min}(g|V):= \min(E(g|V))$ and  $\la_{max}(am):= \max(E(g|V))$ the minimum and maximum.
Define
$$
\la(am) := \frac{\la_{min}(a|\n)}{\la_{max}(m^{-1}|\g)}.
$$
We will construct a function on the set
$$
(AM)^{\sim} \ :=\ \{ am\in AM | \la(am)>1 \}.
$$
The following properties of $(AM)^{\sim}$ are immediate
\begin{itemize}
\item[1.]
$A^+M_{ell}\subset (MA)^{\sim}$
\item[2.]
$am\in (AM)^{\sim} \Rightarrow a\in A^+$
\item[3.]
$am, a'm' \in (AM)^{\sim}, g\in G\ {\rm with}\ a'm'=gamg^{-1} \Rightarrow a=a', g\in AM$.
\end{itemize}
Here we have written $M_{ell}$ for the set of elliptic elements in $M$.

On $(AM)^{\sim}$ define the function $am\mapsto \tilde{l}_{am}$ by
$$
\tilde{l}_{am} \ :=\ \rez{|\alpha|} \log(\la(am)).
$$
Here $\alpha$ is the short positive root in the root system $\Phi(\a ,\g)$.
Note that $\det(1-(am)^{-1}|\n)=0$ implies $\tilde{l}_{am}=0$, since then $a$ and $m^{-1}$ have a common eigenvalue on $\n$.
Note that for $m$ elliptic we have $l_{am}=l_a=|\log(a)|$, the {\bf length} of $a$.
The function $\tilde{l}$ is invariant under conjugation in $AM$ and smooth on a dense open subset set of  $(AM)^{\sim}$.
Let $l_{am}$ be a conjugation invariant smooth function on $(AM)^{\sim}$ such that for any $X$ in the universal enveloping algebra of $\a +\m$ and any $am\in (AM)^{\sim}$ at which $\tilde{l}$ is smooth we have $Xl_{am}\le C_X X\tilde{l}_{am}$, where $C_X$ is a constant only depending on $X$.
Further we insist that for $m$ elliptic we have $\tilde{l}_{am}=l_{am}$.

For $s\in \C$ and $j\in \N$ we define the function $g_s^j$ on $(AM)^{\sim}$ given by
$g_s^j(am) = l_{am}^{j+1} e^{-sl_{am}}$.

Choose any smooth $\eta : N \ra \R$ which has compact support, is positive, invariant under $K\cap M$ and such that
$\int_N \eta(n) dn =1$.

Given these data, let $\Phi =\Phi_{\eta, \tau ,j,s} : G\ra \C$ be defined by
$$
\Phi(kn ma (kn)^{-1}) \= \eta(n) f_\tau(m) \frac{g_s^j(am)}{\det (1-(ma)^{-1}|\n)},
$$
for $k\in K$, $n\in N$, $am\in (AM)^{\sim}$. Further $\phi(g)=0$ if $g$ is not in $G.(MA)^{\sim}$.

In order to see that $\Phi$ is well-defined recall first that by the decomposition $G=KP=KNMA$ every $g\in G.(MA)^{\sim}$ can be written in the form $kn ma (kn)^{-1}$.
By the properties of $(AM)^{\sim}$ we see that two those representations can only differ by an element of $K\cap P =K\cap M$.

The invariance of $\eta$ and of $f_\sigma$ under $K\cap M$-conjugation shows that $\Phi$ is well-defined with respect to the $K$-conjugation.
The points $am\in AM$ where $\det(1-(am)^{-1}|\n)=0$ do not produce poles by the construction of $g_s^j$.

To fix notations let $(\ph ,V_\ph)$ be a finite dimensional unitary representation of $\Ga$ and recall that the group $G$ acts unitarily by right shifts on the Hilbert space $L^2(\Ga \bs G,\ph)$ consisting of all measurable functions $f: G \ra V_\ph$ such that $f(\ga x) = \ph(\ga) f(x)$ and $f$ is square integrable over $\Ga \bs G$ (modulo null functions). 
Let $R$ denote the corresponding representation of $G$.
It is known that this representation splits as
$$
L^2(\Ga \bs G,\ph) \= \bigoplus_{\pi \in \hat{G}} N_{\Ga ,\ph}(\pi) \pi
$$
with finite multiplicities $N_{\Ga ,\ph}(\pi)<\infty$.

Consider a function $f$ on $G$ which is $\dim G +1$-times continuously differentiable and of compact support.
Then $f$ may be plugged into the {\bf Selberg trace formula}:
$$
\tr R(f) \= \sum_{\pi\in\hat{G}} N_{\Ga ,\ph} (\pi)\ \tr \pi(f) \= \sum_{[\ga]}\tr \ph(\ga) \ \vol (\Ga_\ga \bs G_\ga)\ \CO_\ga (f),
$$
where the sum runs over all conjugacy classes in the group $\Ga$ and $\CO_\ga$ is the {\bf orbital integral}:
$$
\CO_\ga(f) \= \int_{G/G_\ga} f(x\ga x^{-1})\ dx.
$$

\begin{proposition} \label{goes_in_tf}
The function $\Phi$ is $(j-\dim(\n))$-times continuously differentiable. For $j$ and $\Re(s)$ large enough it goes into the trace formula  and for any finite dimensional unitary representation $\ph$ of $\Ga$ we have
$$
\sum_{\pi \in \hat{G}} N_{\Ga,\ph} (\pi)\ \tr\pi(\Phi) \= \sum_{[\ga]} \vol (\Ga_\ga \bs G_\ga )\ \O_{m_\ga}^M(f_\tau)\ \frac{g_s^j(a_\ga) \tr (\ph(\ga))}{\det (1-(m_\ga a_\ga)^{-1}|\n)},
$$
where the sum runs over all classes $[\ga]$ such that $\ga$ is conjugate to an element $m_\ga a_\ga$ of $M_{ell}A^+$.
\end{proposition}

\prf
Since the zeroes of $l_{am}$ cancel those of $\det(1-(am)^{-1}|\n)$ it is clear that $\Phi$ is $j-\dim(\n)$-times continuously differentiable.
So we assume $j\ge \dim(\n)$ from now on.
The proof that $\Phi$ goes into the trace formula is nearly wordwise the same as in the rank one case \cite{schub}, \$ 2.

It remains to compute the orbital integrals to see that they give what is claimed in the proposition.
Since $\Phi$ vanishes on all $g\in G$ which are not conjugate to some $am\in A^+M$ it remains to compute the orbital integrals $\CO_{am}(\Phi)$.
The group $\Ga$ being cocompact, only contains semisimple elements therefore it suffices to compute $\CO_{am}(\Phi)$ for semisimple $m$.
We have the decomposition $G=KNAM$ and, in this order, the Haar-measure of $G$ is just the product of the Haar measures of the smaller subgroups.

Since $\Phi(am)=0$ unless $a$  is regular we restrict to $am \in A^+M$, $m$ semisimple.

At first we consider the case, when $m$ is not elliptic and we want to show that the orbital integral $\CO_{am}(\Phi)$ vanishes in this case.
We show more sharply that the integral
$$
\int_{G/AM_m}\Phi(xamx^{-1})\ dx
$$
vanishes whenever $\Phi$ is nontrivial on the orbit of $am$.
Since $AM_m$ is a subgroup of $G_{am}$ this implies the assertion.
The latter integral equals
$$
\int_{KN} \int_{M/M_m} \Phi(kn a m'mm'^{-1}(kn)^{-1})\ dm' dk dn,
$$
which by definition of $\Phi$ equals
$$
\CO_m^M(f_\tau) \frac{g_s^j(am)}{\det(1-(am)^{-1}|\n)}.
$$
Thus we see that $\CO_{am}(\Phi)=0$ unless $m$ is elliptic since this holds for the $M$-orbital integrals over $f_\tau$.

Now we consider $am\in A^+M$ with $m$ elliptic.
In this case we have the equality 
$$
G_{am}\= AM_m.
$$
So that as above we get
$$
\CO_{am}(\Phi) \= \CO_m^M(f_\tau)  \frac{g_s^j(a)}{\det(1-(am)^{-1}|\n)}.
$$
\qed

Besides the parabolic $P=MAN$ we also consider the opposite parabolic $\bar{P} =MA\bar{N}$. 
The complexified Lie algebra of ${N}$ is written ${\n}$. Let $V$ denote a Harish-Chandra module of $G$ then we consider the Lie algebra homology $H_*({\n},V)$ and cohomology $H^*({\n},V)$. It is shown in \cite{HeSch} that these are Harish-Chandra modules of the group $MA$.
For any irreducible unitary representation $\pi$ of $G$ we will denote by $\pi_K$ the Harish-Chandra module of $K$-finite vectors in $\pi$.

We will say that a discrete subgroup $\Ga \subset G$ is {\bf neat} \label{neat} if for every $\ga \in \Ga$ the adjoint $\Ad(\ga)$ has no roots of unity as eigenvalues. 
Every arithmetic $\Ga$ has a neat subgroup of finite index \cite{Bor}.

Let $H_1\in \a^+$ be the unique element with $B(H_1)=1$.

We will denote by ${\cal E}_P(\Ga)$ the set of nontrivial 
$\Ga$-conjugacy classes $[\ga]$, which are such that $\ga$ is 
in $G$ conjugate to an element of $A^+B$.
Such an element will then 
be written $a_\ga b_\ga$ or $a_\ga m_\ga$ with $a_\ga \in A^+$. 
The element $\ga \neq 1$ will be called {\bf primitive} if 
$\sigma \in \Ga$ and $\sigma^n =\ga$ with $n\in \N$ implies $n=1$. 
Every $\ga \neq 1$ is a power of a unique primitive element. 
Obviously primitivity is a property of conjugacy classes. 
Let $\CE_P^p(\Ga)$ denote the subset of $\CE_P(\Ga)$ consisting 
of all primitive classes.

The following generalizes results of A. Juhl \cite{Ju}.

\begin{theorem}\label{genSelberg}
Let $\Ga$ be neat and $(\ph ,V_\ph)$ a finite dimensional unitary 
representation of $\Ga$. Choose a parabolic $P=MAN$ of 
splitrank one. For $\Re(s)>>0$ define the {\bf generalized Selberg zeta function}:
$$
Z_{P,\tau,\ph}(s) \= \prod_{[\ga]\in {\cal E}_P^p(\Ga)} \prod_{N\geq 0} 
\det\left(1-e^{-sl_\ga}\ga \left| \begin{array}{c}V_\ph \otimes V_\tau 
\otimes S^N(\n)\end{array}\right.\right)^{\chi_{_1}(X_\ga)},
$$
where $S^N(\n)$ denotes the $N$-th symmetric power of the space $\n$ and $\ga$ acts on 
$V_\ph \otimes V_\tau \otimes S^N(\n)$ via $\ph(\ga) \otimes \tau(m_\ga) 
\otimes Ad^N((m_\ga a_\ga)^{-1})$, here $\ga \in \Ga$ is conjugate to 
$\m_\ga a_\ga \in BA^+$. 

Then $Z_{P,\tau,\ph}$ has a meromorphic 
continuation to the entire plane. 
The vanishing order of $Z_{P,\tau ,\ph}(s)$ at a point $s=\la (H_1)$, $\la \in \a^*$, is
$$
(-1)^{\dim \ \n} \sum_{\pi \in \hat{G}}N_{\Ga ,\ph}(^\theta\pi)
\sum_{p,q}(-1)^{p+q} \dim \left(\begin{array}{c}H^q({\n},\pi_K)\otimes \wedge^p\p_M \otimes V_{\tau}\end{array}\right)^{K_M}_{-\la},
$$
where $(.)_\la$ denotes the generalized $\la$-eigenspace.
Here $\theta$ is the Cartan involution which 
acts on $\hat{G}$ by $^\theta\pi(g)=\pi(\theta(g))$.
It follows that all poles and zeroes of the function $Z_{P,\tau ,\ph}(s+|\rho_0|)$ lie in $\R \cup i\R$.
\end{theorem}

$ $

\prf
From \cite{D-Hitors} we take for $\Ga$ neat
$$
\chi_{_1}(X_\ga) \= \frac{| W(\g_\ga ,\h)| \prod_{\alpha \in \Phi_\ga^+}(\rho_\ga ,\alpha )}{l_{\ga_0} c_\ga [G_\ga :G_\ga^0]} \vol (\Ga_\ga \backslash G_\ga),
$$
where $l_{\ga_0}$ is the length of the primitive $\ga_0$ underlying $\ga$. 

Now assume $\ga\in\CE_P(\Ga)$, so $\ga$ is conjugate to $am\in A^+B$.
Then $G_\ga$ is conjugate to $A M_{m}$.
Since $A$ is connected and meets $M$ only in the trivial element it follows $[G_\ga :G_\ga^0] =[M_m:M_m^0]$ and with Proposition \ref{orbitalint} we conclude
$$
\CO_m(f_\tau)\vol(\Ga_\ga\bs G_\ga) = l_{\ga_0} \tr \tau(m) \chi_1(X_\ga).
$$

Since $\Ga$ is neat and $G$ is of adjoint type it follows that for any $\ga\in\Ga$ we have for the centralizers: $G_{\ga^n}=G_\ga$ for any natural number $n$.
In \cite{DKV} it is proven that $X_\ga \cong \Ga_\ga \bs G_\ga /K_\ga$ and so it follows that $X_{\ga^n} \cong X_{\ga}$ for any $n$.
With this in mind we get from Proposition \ref{goes_in_tf} that the right hand side of the trace formula is:
$$
\sum_{[\ga]} \frac{l_{\ga_0} \chi_{_1}(X_\ga) {\tr\ \tau(m_\ga)} l_\ga^{j+1} e^{-sl_\ga}\tr(\ph(\ga))}{\det(1-\ga^{-1} | \n)}
\= (-1)^{j+1} (\der{s})^{j+2} \log Z_{P,\tau,\ph} (s).
$$

To understand the spectral side let
$$
(MA)^+ \= {\rm interior\ in}\ MA\ {\rm of\ the\ set}\ \{ g\in MA | \det(1-ga|\bar{\n})\geq 0\ \forall a\in A^+ \}.
$$
Note that, if the real rank of $G$ is $1$, then $(MA)^+$ coincides with the set of  all $g\in MA$ such that $\det (1-g|\bar{\n})>0$ which also equals $MA^+$.
This becomes false in the higher rank case.
Note also that always $M_{ell}A^+\subset (MA)^+$ holds.
Now $B$ is a compact Cartan in $M$ and $(M.B)A^+$ is a subset of $(MA)^+$. 
Let $\pi\in\hat{G}$ then the space of $K$-finite vectors $\pi_K$ forms a Harish-Chandra-module.
Any Harish-Chandra-module $V$ has a character which by the work of Harish-Chandra is known to be a function $\Theta_V^G$. 
In the case of $\pi_K$ we write that function as $\Theta_\pi^G$.
It was proven in \cite{HeSch} that
 for any $\pi \in \hat{G}$ we have for $ma \in (MA)^+ \cap G^{\rm reg}$
$$
\Theta_\pi^G(ma) \= \frac{\Theta_{H_*(\bar{\n},\pi_K)}^{MA}(ma)}{\det(1-ma|\bar{\n})}.
$$
Let $f$ be integrable and supported on $G.(MA)$ then applying the Weyl integration formula first to $G$ and then backwards to $M$ gives
$$
\int_G f(g)\ dg \= \int_{G/MA} \int_{MA^+} f({g}ma{g}^{-1})\ |det(1-ma|\n \oplus \bar{\n})|\ da dm d{g}.
$$
So that for $\pi \in \hat{G}$
\begin{eqnarray*}
\tr\ \pi(\Phi) &\=& \int_G \Theta_\pi^G(x) \Phi(x)\ dx
\\
        &\=& \int_{MA^+} \Theta_\pi^G(ma)\ f_\tau(m)\ 
 \left| \frac{\det(1-ma|{\n}\oplus \bar{\n})}{\det(1-(ma)^{-1}|{\n})}\right| \ g_s^j(am)\ da dm.
\end{eqnarray*}
Applying the Weyl integration formula again we see that the integral over $M$ is a sum over all conjugacy classes of Cartan subgroups $L$ of $M$ where the summand is $\rez{|W(L,M)|}$ times
$$
\int_L\int_{M/L} f_\tau({m}l{m}^{-1})\ d{m}\  D(l)\ \Theta_\pi^G(al)\ \left|\frac{\det(1-al|{\n}\oplus \bar{\n})}{\det(1-(al)^{-1}|{\n})}\right|\ g_s^j(al)\ dl,
$$
here $D(l)$ is the Jacobian determinant of the Weyl integration formula on $M$.
The integral over $M/L$ is just an orbital integral of $f_\tau$ so that in the sum only the compact Cartan $B$ will survive. Hence we get that
\begin{eqnarray*}
\tr \pi (\Phi) &\=& \rez{|W(B,M)|} \int_{A^+B} 
	\int_{M/B} f_\tau(mbm^{-1})\ dm \  \Theta_{H_*(\bar{\n},\pi_K)}^{MA} (ab)
\\
{}&{}&\ \ \ \ \ \ \times \left| \frac{\det(1-ab |\n)}
			{\det(1-(ab)^{-1} |\n)}\right|\ g_s^j (ab)\ da db
\end{eqnarray*}
By the compactness of $B$ the Jacobian factor here is trivial. The factor 
$$
\left| \frac{\det(1-ab |\n)}
	    {\det(1-(ab)^{-1} |\n)}\right|
$$
equals $\det(a|\wedge^{max}\n)$.
Using the isomorphism $H_p(\bar{\n} ,V) \cong H^{\dim \ \bar{\n} -p}(\bar{\n} ,V) \otimes \wedge^{max} \bar{\n}$
we replace 
$$
\Theta_{H_*(\bar{\n},\pi_K)}^{MA} (ab)\det(a|\wedge^{max}\n)
$$
by
$$
(-1)^{\dim \n}	\Theta_{H^*(\bar{\n},\pi_K)}^{MA} (ab).
$$
A final application of the Weyl integration formula yields
$$
\tr\ \pi(\Phi) \= (-1)^{\dim \ \n} \int_{MA^+} f_{{\tau}}(m)\ \Theta_{H^*(\bar{\n}, \pi_K)}^{MA}(ma)\ dm g_s^j(a)\ da.
$$
Using the properties of $f_\tau$ we get
$$
\tr\ \pi (\Phi) \= (-1)^{\dim \ \n}\int_{A^+} \tr (a|(H^*(\bar{\n},\pi_K) \otimes \wedge^* \p_M \otimes V_{\tau})^{K_M})\ l_a^{j+1} e^{-sl_a} \ da
$$
Now let the $A$-decomposition of $(H^*(\bar{\n},\pi_K) 
\otimes \wedge^* \p_M \otimes V_{\tau})^{K_M}$ be 
$\bigoplus_{\la\in\La_\pi} E_\la,$ where $E_\la$ is the biggest subspace on 
which $a-e^{\la(\log a)}$ acts nilpotently. 

Let $m_\la$ denote the 
dimension of $E_\la$.
Integrating over $\a^+$ instead over $A^+$ we get
\begin{eqnarray*}
\tr\ \pi (\Phi) &\=& (-1)^{\dim (\n)}\int_0^\infty \sum_\la m_\la e^{(\la(H_1)-s)t} t^{j+1}\ dt
\\
        &\=& (-1)^{\dim (\n) +j+1} (\der{s})^{j+1} \sum_\la m_\la \rez{s-\la(H_1)}
\end{eqnarray*}

From this it follows that the order of $Z_{P,\tau ,\ph}(s)$
at a point $s=\la (H_1)$, $\la \in \a^*$, is
$$
(-1)^{\dim \ \n} \sum_{\pi \in \hat{G}}N_{\Ga ,\ph}(\pi)
\sum_{p,q}(-1)^{p+q} \dim \left(\begin{array}{c}H^q(\bar{\n},\pi_K)
\otimes \wedge^p\p_M \otimes V_{\tau}\end{array}\right)^{K_M}_{\la}.
$$

The Cartan-involution $\theta$ maps $\bar{\n}$ to $\n$.
Applying this and summing over $^\theta\pi$ we get
$$
(-1)^{\dim \ \n} \sum_{\pi \in \hat{G}}N_{\Ga ,\ph}(^\theta\pi)
\sum_{p,q}(-1)^{p+q} \dim \left(\begin{array}{c}H^q({\n},\pi_K)\otimes 
\wedge^p\p_M \otimes V_{\tau}\end{array}\right)^{K_M}_{-\la},
$$
and we are done.

To prove the generalized Riemann hypothesis for this zeta function recall that in the tensor product 
$$
H^q({\n},\pi_K)\otimes \wedge^p\p_M \otimes V_{\tau}
$$ 
the 
group $A$ only acts on the first tensor factor.
By Corollary 3.32 in \cite{HeSch} we know that the possible eigenvalues 
of $A$ on $H^q({\n},\pi_K)$ are given by $w\la_\pi +\rho|_\a$, where $\la_\pi$
is the infinitesimal character of the representation $\pi$ and $w$ is any 
element of the Weyl group $W(\g ,\h)$.
For any $\pi \in \hat{G}$ that contributes, the global character $\Theta_\pi$ must be nonvanishing on $G.H^{reg}$.
Any global character is a linear combination of global characters induced from tempered representations (\cite{Knapp} p.383.)
Since $P$ is a maximal parabolic we see that $\pi$ must be in the discrete series, the limits of the discrete series or induced 
from the discrete or limit series of $M$. In any case the infinitesimal character $\la_\pi$ is real on $\b$.
Now $\pi (C)$ is selfadjoint, where $C$ is the Casimir operator, and by the well known $\pi (C) =(B(\la_\pi)-B(\rho))Id$ we get that $B(\la_\pi)$ is real. But this is $B(\la_\pi |_\a)+B(\la_\pi |_\b)$ and the latter is real so $B(\la_\pi |_\a)$ is real.
Since $\a$ is one dimensional, $\la_\pi |_\a$ must be real or purely imaginary.
\qed

In case that the Weyl group $W(G,A)$ is nontrivial, we have $Z_{P,\tau,\ph}=Z_{\theta(P),\tau,\ph}$.
To give the reader a feeling of this condition consider the case $G=SL_3(\R)$. 
In that case we have trivial Weyl group $W(G,A)$.
On the other hand, consider the case when the fundamental rank of $G$ is $0$; this is the most interesting case to us since only then we have several conjugacy classes of splitrank-one Cartan subgroups.
In that case it follows that the dimension of all irreducible factors of the symmetric space $X=G/K$ is even, hence the point-reflection at the point $eK$ is in the connected component of the group of isometries of $X$. 
This reflection can be thought of as an element of $K$ which induces a nontrivial element of the Weyl group $W(G,A)$.
So we see that in the most important case we have $|W(G,A)|=2$.

From now on for the rest of the section we assume that the Weyl group $W(G,A)$ is nontrivial.
Let $w$ denote the nontrivial element as well as a representative in $K$.
We will restrict to the case when the representation $\tau$ is invariant under $w$. This means $\tau(m)=^w\tau(m)=\tau(wmw^{-1})$.
Since $\tau$ need not be irreducible this can be achieved by replacing $\tau$ by $\tau \oplus ^w\tau$ if necessary.
For $p\in \{ 0,\dots ,\dim (\p_M)+\dim (N)\}$ let $V_p$ be the space
$$
V_p \= \bigoplus_{a+b=p} \left( H^a({\n},C^\infty(\Ga \bs G,\ph )_{K-{\rm finite}})
\otimes \wedge^b \p_M \otimes V_{\tau}\right)^{K_M}.
$$
The normalized element $H=H_1$ of $\a$ acts on $V_p$.

\begin{proposition} \label{detformel}
Assume $\tau$ is invariant under $w$, then
the zeta function satisfies a determinant formula, i.e.:
$$
Z_{P,\tau ,\ph}(s) \= e^{Q(s)} \det 
\left( H+s\left|\bigoplus_p (-1)^pV_p\right.\right)^{(-1)^{\dim (N)}} ,
$$
where $Q$ is a polynomial of degree $\leq \dim G +\dim N$.
\end{proposition}

\prf
Suppose we could show the existence of the determinant on the right hand side.
The proof of the theorem then shows that the derivatives of order 
$>\dim G +\dim N$ of the logarithm of the RHS coincide with the 
corresponding terms of the left hand side. The assertion follows.

We are left to show the existence of the determinant.
For this we will employ Cram\'er's theorem \cite{JoLa}, which says 
that $Z_{P,\tau ,\ph}$ is a determinant if it satisfies a functional 
equation with fudge factor being a determinant.

A representation $\pi\in\hat{G}$ is called {\bf elliptic} if $\Theta_\pi$ is 
nonzero on the compact Cartan. 
Let $\hat{G}_{ell}$ be the set of elliptic elements in $\hat{G}$ and 
denote by $\hat{G}_{ds}$ the subset of discrete series representations.
Further let $\hat{G}_{lds}$ denote the set of all discrete series and all limits of discrete series representations.
In the proof of the theorem we have shown that the vanishing order of 
$Z_{P,\tau ,\ph}(s)$ at the point $s=\mu (H_1)$, $\mu\in\a^*$ is
$$
(-1)^{\dim N}\sum_{\pi\in\hat{G}} N_{\Ga ,\ph}(\pi) m(\pi ,\tau ,\mu),
$$
where
$$
m(\pi ,\tau ,\mu) \= \sum_{p,q} (-1)^{p+q} \dim 
\left(\begin{array}{c}H^q(\n ,\pi_K)\otimes \wedge^p\p_M\otimes 
V_{\tau}\end{array}\right)_{-\mu}^{K_M}.
$$

A standard representation is a representation induced from a tempered representation.
Any character $\Theta_\pi$ for $\pi\in\hat{G}$ is an integer linear combination
of characters of standard representations.
From this it follows that for $\pi\in\hat{G}$ the character restricted to 
the compact Cartan $T$ is
$$
\Theta_\pi \mid_T \= \sum_{\pi'\in\hat{G}_{lds}} k_{\pi ,\pi'} \Theta_{\pi'}
\mid_T,
$$
with integer coefficients $k_{\pi ,\pi'}$.

\begin{lemma}
Assume $\tau$ is invariant under $w$, then there is a $C>0$ such that for $\Re (\mu(H_1))<-C$ the order of 
$Z_{P,\tau ,\ph}(s)$ at $s=\mu(H_1)$ is
$$
(-1)^{\dim N}\sum_{\pi\in\hat{G}_{ell}} N_{\Ga ,\ph}(\pi) 
\sum_{\pi'\in\hat{G}_{lds}} k_{\pi ,\pi'} m(\pi' ,\tau ,\mu).
$$
\end{lemma}

\prf
For any $\pi\in\hat{G}$ we know that if $\Theta_\pi |_{AT_M}\neq 0$ 
then in the representation of $\Theta_\pi$ as linear combination of
standard characters there must occur lds-characters and characters
of representations $\pi_{\xi ,\nu}$ induced from $P$.
Since $\Theta_{\pi_{\xi ,\nu}}=\Theta_{\pi_{^w\xi ,-\nu}}$ and $\tau$ is $w$-invariant any contribution 
of $\pi_{\xi ,\nu}$ for $\Re (s) <<0$ would also give a pole or zero of 
$Z_{P,\tau ,\ph}$ for $\Re(s)>>0$.
In the latter region we do have an Euler product, hence there are 
no poles or zeroes.
\qed

Now consider the case $FR(G)=1$, so there is no compact Cartan, hence no discrete series.
Since $\tau$ is invariant under $w$ it follows that the poles and zeroes of $Z_{P,\tau ,\ph}(s+| \rho_0|)$ 
lie symmetric around $s=0$.
By a finite genus argument it follows that

\begin{proposition} There is a polynomial $P$ of degree $\leq \dim G+\dim N$ such that
$$
Z_{P,\tau ,\ph}(s) \= e^{P(s)} Z_{P,\tau ,\ph}(2| \rho_0| -s).
$$
\end{proposition}
\qed

The claim of Proposition \ref{detformel} follows in this case.

Now assume $FR(G)=0$ so there is a compact Cartan subgroup $T$.
As Haar measure on $G$ we take the Euler-Poincar\'e measure.
The sum in the lemma can be rearranged to
$$
(-1)^{\dim N}
\sum_{\pi'\in\hat{G}_{lds}}  m(\pi' ,\tau ,\mu)\sum_{\pi\in\hat{G}_{ell}} 
N_{\Ga ,\ph}(\pi) k_{\pi ,\pi'}.
$$
We want to show that the summands with $\pi'$ in the limit of the discrete
series add up to zero.
For this suppose $\pi'$ and $\pi''$ are distinct and belong to the limit of the discrete series. Assume further that their Harish-Chandra parameters agree.
By the Paley-Wiener theorem \cite{CloDel} there is a smooth compactly supported function $f_{\pi' ,\pi''}$ such that for any tempered $\pi\in\hat{G}$:
$$
\tr \pi(f_{\pi' ,\pi''}) \= \left\{ \begin{array}{cl} 1& {\rm if}\ \pi =\pi'\\
						-1 	&{\rm if}\ \pi =\pi''\\
						0	& {\rm else.}
					\end{array}\right.
$$
Plugging $f_{\pi' ,\pi''}$ into the trace formula it follows
$$
\sum_{\pi\in\hat{G}_{ell}} N_{\Ga ,\ph}(\pi)k_{\pi ,\pi'}
\=\sum_{\pi\in\hat{G}_{ell}} N_{\Ga ,\ph}(\pi)k_{\pi ,\pi''},
$$
so that in the above sum the summands to $\pi'$ and $\pi''$ occur with the same coefficient.
Let $\pi_0$ be the induced representation whose character is the sum of the characters of the $\pi''$, where $\pi''$ varies over all lds-representations with the same Harish-Chandra parameter as $\pi'$.
Then for $\Re (\mu(H_1))<-C$ we have $m(\pi_0,\tau,\mu)=0$.
Thus it follows that the contribution of the limit series vanishes.	

Plugging the pseudo-coefficients \cite{Lab} of the discrete series 
representations into the trace formula gives for $\pi\in\hat{G}_{ds}$:
$$
\sum_{\pi'\in\hat{G}_{ell}} k_{\pi' ,\pi} N_{\Ga ,\ph}(\pi') 
\= \dim \ph (-1)^{\frac{\dim X}{2}}\chi(X_\Ga) d_{\pi},
$$
where $d_{\pi}$ is the formal degree of $\pi$.

The infinitesimal character $\la$ of $\pi$ can be viewed as an element of the coset space $(\t^*)^{reg}/W_K$.
So let $J$ denote the finite set of connected components of $(\t^*)^{reg}/W_K$, then we get a decomposition $\hat{G}_{ds} = \coprod_{j\in J} \hat{G}_{ds,j}$.

In the proof of the theorem we used the Hecht-Schmid character 
formula to deal with the global characters.
On the other hand it is known that global characters are given on the regular
set by sums of toric characters over the Weyl denominator.
So on $H=AB$ the character $\Theta_\pi$ for $\pi\in\hat{G}$ is of the form
$\CN /D$, where $D$ is the Weyl denominator and the numerator $\CN$ is of the form
$$
\CN (h) \= \sum_{w\in W(\t,\g)} c_w h^{w\la},
$$
where $\la\in\h^*$ is the infinitesimal character of $\pi$.
Accordingly, the expression $m(\pi ,\tau ,\mu)$ expands as a sum 
$$
m(\pi ,\tau ,\mu) \= \sum_{w\in W(\t ,\g)} m_w(\pi ,\tau ,\mu).
$$

\begin{lemma}
Let $\pi ,\pi'\in \hat{G}_{ds,j}$ with infinitesimal characters $\la ,\la'$ 
which we now also view as elements of $(\h^*)^+$, then
$$
m_w(\pi ,\tau ,\mu) \= m_w(\pi',\tau ,\mu +w(\la'-\la)|_\a).
$$
\end{lemma}

\prf
In light of the preceding it suffices to show the following:
Let $\tau_\la, \tau_{\la'}$ denote the numerators of the global characters 
of $\pi$ and $\pi'$ on $\h^+$.
Write
$$
\tau_\la (h) \= \sum_{w\in W(\t ,\g)}c_w h^{w\la}
$$
for some constants $c_w$. Then we have
$$
\tau_{\la'} (h) \= \sum_{w\in W(\t ,\g)}c_w h^{w\la'}.
$$
To see this, choose a $\la''$ dominating both $\la$ and $\la'$, then apply the Zuckerman functors $\ph_{\la''}^\la$ and $\ph_{\la''}^{\la'}$.
Proposition 10.44 of \cite{Knapp} gives the claim.
\qed

Write $\pi_\la$ for the discrete series representation with infinitesimal 
character $\la$.
Let $d(\la):= d_{\pi_\la}$ be the formal degree then $d(\la)$ is a polynomial 
in $\la$, more precisely from \cite{AtSch} we take
$$
d(\la) \= \prod_{\alpha\in\Phi^+(\t ,\g)} 
\frac{(\alpha ,\la +\rho)}{(\alpha ,\rho)},
$$
where the ordering $\Phi^+$ is chosen to make $\la$ positive.

Putting things together we see that for $\Re(s)$ small enough the order of 
$Z_{P,\tau ,\ph}(s)$ at $s=\mu(H_1)$ is
\begin{eqnarray*}
\CO (\mu) &\=& \dim \ph (-1)^{\frac{\dim X}{2}} \chi (X_\Ga)\\
	& {}&\sum_{j\in J} \sum_{w\in W(\t ,\g)} \sum_{\pi \in \hat{G}_{ds,j}}
	d(\la_\pi) m_w(\pi_j,\tau ,\mu +w(\la_j-\la_\pi)|_\a),
\end{eqnarray*}
where $\pi_j\in\hat{G}_{ds,j}$ is a fixed element.
The function $\mu \mapsto m_w(\pi_j ,\tau ,\mu)$ takes nonzero values only for 
finitely many $\mu$.
Since further $\la\mapsto d(\la)$ is a polynomial it follows that the 
regularized product
$$
D_{H,\tau ,\ph}(s) := \hat{\prod_{\mu ,\CO(\mu)\neq 0}} (s-\mu(H_1))^{\CO(\mu)}
$$
exists. Further it is clear that with
$$
\hat{Z}_{H,\tau ,\ph}(s) := Z_{P,\tau ,\ph}(s) D_{H,\tau ,\ph}(s)^{-1}
$$
we have
$$
\hat{Z}_{H,\tau ,\ph}(2| \rho_0| -s) \= e^{P(s)} \hat{Z}_{H,\tau ,\ph}(s),
$$ where $P$ is a polynomial.
By the theorem of Cram\'er Proposition \ref{detformel} follows.

\qed

\begin{proposition} \label{extformel}
Let $\sigma$ be a finite dimensional representation of $M$ then the order of
$Z_{P,\sigma , \ph}(s)$ at $s=\la(H_1)$ is
$$
(-1)^{\dim(N)} \sum_{\pi \in \hat{G}}N_{\Ga,\ph}(^\theta\pi) \sum_{q=0}^{\dim(\m\oplus{\n} /\k_M)}(-1)^q \dim(H^q(\m\oplus{\n}, K_M ,\pi_K \otimes V_{{\sigma}})_{-\la}).
$$
This can also be expressed as
$$
(-1)^{\dim(N)} \sum_{\pi \in \hat{G}}N_{\Ga,\ph}(^\theta\pi) \sum_{q=0}^{\dim(\m\oplus{\n} /\k_M)}(-1)^q \dim({\rm Ext}_{(\m \oplus {\n} ,K_M)}^q (V_{\breve{\sigma}} ,V_\pi)_{-\la}).
$$
\end{proposition}

\prf
Extend $V_{{\sigma}}$ to a $\m \oplus {\n}$-module by letting ${\n}$ act trivially. We then get
$$
H^p({\n},\pi_K) \otimes V_{\breve{\sigma}} \cong H^p({\n},\pi_K \otimes V_{\sigma}).
$$

The $(\m ,K_M)$-cohomology of the module $H^p({\n},\pi_K \otimes V_{\sigma})$ is the cohomology of the complex $(C^*)$ with
\begin{eqnarray*}
C^q &\=& {\rm Hom}_{K_M}(\wedge^q\p_M ,H^p({\n},\pi_K)\otimes V_{\sigma})
\\
        &\=& (\wedge^q\p_M \otimes H^p({\n},\pi_K)\otimes V_{\sigma})^{K_M},
\end{eqnarray*}
since $\wedge^p\p_M$ is a self-dual $K_M$-module. Therefore we have an isomorphism of virtual $A$-modules:
$$
\sum_q (-1)^q (H^p({\n},\pi_K)\otimes \wedge^q\p_M \otimes V_{\sigma} )^{K_M}
\cong
\sum_q (-1)^q H^q (\m ,K_M,H^p({\n},\pi_K \otimes V_{\sigma})).
$$

Now one considers the Hochschild-Serre spectral sequence in the relative case for the exact sequence of Lie algebras
$$
0 \ra {\n}\ra \m \oplus {\n}\ra \m \ra 0
$$
and the $(\m\oplus {\n},K_M)$-module $\pi \otimes V_{\sigma}$. We have
$$
E_2^{p,q} \= H^q(\m ,K_M ,H^p({\n},\pi_K\otimes V_{\sigma}))
$$
and
$$
E_\infty^{p,q} \= {\rm Gr}^q(H^{p+q}(\m \oplus {\n},K_M ,\pi_K \otimes V_{\sigma})).
$$
Now the module in question is just
$$
\chi(E_2) \= \sum_{p,q} (-1)^{p+q} E_2^{p,q}.
$$
Since the differentials in the spectral sequence are $A$-homomorphisms this equals $\chi(E_\infty)$.
So we get an $A$-module isomorphism of virtual $A$-modules
$$
\sum_{p,q} (-1)^{p+q} (H^p({\n},\pi_K)\otimes \wedge^q\p_M \otimes V_{\sigma})^{K_M} \cong \sum_j (-1)^j H^j (\m\oplus {\n},K_M, \pi_K \otimes V_{\sigma}).
$$
The second statement is clear by \cite{BorWall} p.16. \qed

Consider the space
$$
W_p \= H^p(\m \oplus {\n},K_M,C^\infty(\Ga \bs G,\ph)_{K-{\rm finite}}\otimes V_{\sigma}).
$$
The space $W_p$ can be thought of as a cohomology space of the tangential cohomology of the unstable-central bundle attached to the flow $\phi_t$ on $\Ga \bs G /K_M$ with generator $H$.
Since the closed geodesics occurring in the zeta function are just the closed orbits of $\phi_t$, the following proposition gives the analogue of the determinant formula in Deligne's proof of the Weil conjectures. In the rank one case the following is already implicitly in \cite{Ju}.

\begin{proposition} 
For the zeta function $Z_{P,\sigma ,\ph}$ attached to a $w$-invariant $M$-representation $\sigma$ we have
$$
Z_{P,\sigma ,\ph}(s) \= e^{Q(s)}
\det \left( H+s\left|\bigoplus_p(-1)^pW_p\right.\right)^{(-1)^{\dim (N)}},
$$
for a polynomial $Q$ of degree $\leq \dim G+\dim N$.
\end{proposition}
\qed

Note that if $W(G,A)$ is nontrivial then $\ga\in\CE_P(\Ga)$ implies $\ga^{-1}\in\CE_P(\Ga)$.
Considering the definition of the zeta function we then find that $Z_{P,\tau ,\ph}=Z_{\theta(P),\tau ,\ph}$, but this does not hold in general.
To get rid of the dependence on $P$ we define the {\bf zeta function attached to the Levi-component $MA$} by
$$
Z_{MA,\tau ,\ph}(s) := Z_{P,\tau ,\ph}(s) Z_{\theta(P),\tau ,\ph}(s).
$$

\section{The Ruelle zeta function}
The Ruelle zeta function can be described in terms of the Selberg zeta function as follows.

\begin{theorem}
Let $\Ga$ be neat and choose a parabolic $P$ of splitrank one. For $\Re(s)>>0$ define the zeta function
$$
Z_{P,\ph}^R(s) \= \prod_{[\ga]\in {\cal E}_H^p(\Ga)} \det\left(\begin{array}{c}1-e^{-sl_\ga}\ph(\ga)\end{array}\right)^{\chi_{_1}(X_\ga)},
$$
then $Z_{P,\ph}^R(s)$ extends to a meromorphic function on $\C$. In the case that there is only one positive root $\alpha$ in the root system of $(\a ,\g)$ we have 
$$
Z_{P,\ph}^R(s) \= \prod_{l=0}^{\dim N} Z_{P,\sigma_l,\ph}(s+l|\alpha |)^{(-1)^l},
$$
where $\sigma_l$ is the representation of $K_M$ on $\wedge^l \n$.
If we have two positive roots, say $\alpha$ and $2{\alpha}$, then we have ${\n} = {\n}_r \oplus {\n}_I$ where ${\n}_r$ has dimension one and $a$ acts on ${\n}_r$ by $a^{2\alpha}$ and on ${\n}_I$ by $a^{2{\alpha}}$.
Then it follows
$$
Z_{P,\ph}^R(s) \= \prod_{l=0}^{\dim N-1} \left( \frac{Z_{P,\wedge^l{\n}_I,\ph}(s+l|\alpha |)}{Z_{P,\wedge^l{\n}_I\otimes {\n}_r,\ph}(s+(l+2)|\alpha |)}\right)^{(-1)^l}. 
$$
\end{theorem}

\prf
Clear.
\qed

\section{$\Ga$-cohomology and the divisor of the zeta function}

In this section we consider a finite dimensional representation $(\sigma , V_\sigma)$ of $M$.

Let $\nu \in \a^*$ and consider the principal series representation $\pi_{\sigma ,\nu}$ on the Hilbert space $H_{\sigma ,\nu}$ of all functions $f$ from $G$ to $V_\sigma$ such that $f(xman)=a^{-(\nu +\rho)} \sigma(m)^{-1} f(x)$ and such that the restriction of $f$ to $K$ is an $L^2$-function.
For $\nu \in \a^*$ let $\bar{\nu}$ denote its complex conjugate with respect to the real form $\a_0^*$.

According to \cite{KaSch} there is to each Harish-Chandra module $V$ a minimal globalization $V^{min}$ and a maximal globalization $V^{max}$.
We will now formulate the adapted Patterson conjecture \cite{BuOl}.

\begin{theorem}
The cohomology $H^p(\Ga ,H_{\sigma ,\nu}^{max} \otimes V_\ph)$ is finite dimensional for all $p\geq 0$ and the vanishing order of the zeta function
$Z_{P,\sigma ,\ph}$ is
\begin{eqnarray*}
\ord_{s=\nu(H_1)}Z_{P,\sigma ,\ph}(s+|\rho_0|)
&\=& \chi_{_1}(\Ga ,H_{\sigma ,{-\nu}}^{max} \otimes V_\ph)\\ 
&\=& - \sum_{p=0}^\infty p(-1)^p
\dim\ H^p(\Ga ,H_{\sigma ,{-\nu}}^{max} \otimes V_\ph)
\end{eqnarray*}
if $\nu \neq 0$ and
\begin{eqnarray*}
\ord_{s=0}Z_{P,\sigma ,\ph}(s+|\rho_0|) 
&\=& \chi_{_1}(\Ga ,\hat{H}_{\sigma,0}^{max}\otimes V_\ph)\\
&\=& -\sum_{p=0}^\infty p(-1)^p \dim\
H^p(\Ga ,\hat{H}_{\sigma ,0}^{max}\otimes V_\ph),
\end{eqnarray*}
where $\hat{H}_{\sigma ,0}^{max}$ is a certain nontrivial extension of $H_{\sigma ,0}^{max}$ with itself.

Further $\chi (\Ga ,H_{\sigma ,\nu}^{max}\otimes V_\ph)=\sum_p(-1)^p\dim\ H^p(\Ga ,H_{\sigma ,\nu}^{max} \otimes V_\ph)$ always vanishes.
\end{theorem}

Before proving this theorem we are going to sketch some consequences. 
For this consider the special case ${\rm rank} (X)=1$ and $\dim(X)$ even.
In \cite{D-det} the author showed that all poles and zeroes of $Z_{P,\sigma ,\ph}$ lie in the union of the set $\{-1,-2,\dots \}$ and the set $[-a,a] \cup i\R$ for some $a>0$ and that for $s_0=\nu(H_1) \in i\R - \{ 0\}$ it holds
$\ord_{s=s_0}Z_{P,\sigma ,\ph}(s+|\rho_0|) = N_{\Ga ,\ph}(\pi_{\sigma ,\nu}).$

On the other hand, the duality theorem (\cite{Gelf} p.61) implies 
$$
N_{\Ga ,\ph}(\pi) \= \dim\ H^0(\Ga,V_\pi^{max}\otimes V_\ph)
$$
for any $\pi \in \hat{G}$. Therefore we also conclude the relation
$$
\sum_{p=1}^{\dim(X)} (p-1)(-1)^p \dim\ H^p(\Ga ,H_{\sigma ,\nu}^{max}\otimes V_\ph) \= 0.
$$

For example in the case $\dim(X)=2$, so $X=$ hyperbolic plane, this implies for $\nu$ imaginary and $\neq 0$ that
$\dim\ H^2(\Ga ,H_{\sigma ,\nu}^{max}\otimes V_\ph)=0$
and
$N_{\Ga ,\ph}(H_{\sigma ,\nu}) = \dim\ H^0(\Ga ,H_{\sigma ,\nu}^{max}\otimes V_\ph) = \dim\ H^1(\Ga ,H_{\sigma ,\nu}^{max}\otimes V_\ph).$
  
{\bf Proof of the theorem:}
In \cite{KaSch} it is proven that $V^{max}$ coincides with $V^{-\omega}$, the hyperfunction globalization.
We will from now on use this globalization.

We will extend $V_\sigma$ to a $MAN$-module $V_{\sigma_\nu}$ by
$\sigma_\nu (man) = a^{\nu +\rho} \sigma(m)$. 
Let $C^{-\omega}(G)$ denote the hyperfunctions on $G$ and consider $C^{-\omega}(G) \otimes V_{\sigma_\nu}$ as a
$P=MAN$-module, where $P$ acts on $C^{-\omega}(G)$ by right translations.
Further write $C^\infty(G), C^{-\infty}(G)$ for the spaces of $C^\infty$-functions and distributions on $G$.

\begin{lemma}
For $\alpha = \infty ,-\infty ,-\omega$ the $(\m \oplus \a \oplus \n ,K_M)$-cohomology complex of $C^\alpha(G)\otimes V_{\sigma_\nu}$:
\begin{eqnarray} \label{complex}
0 \ra (C^\alpha (G) \otimes V_{\sigma_\nu})^{K_M} \ra 
(C^\alpha (G) \otimes V_{\sigma_\nu} \otimes \wedge^1(\p_M \oplus \a \oplus \n)^*)^{K_M}
\end{eqnarray}
\begin{eqnarray*}
\ra (C^\alpha (G) \otimes V_{\sigma_\nu} \otimes \wedge^2(\p_M \oplus \a \oplus \n)^*)^{K_M} \ra \dots
\end{eqnarray*}
is exact in all degrees $p\geq 1$.
\end{lemma}

\prf
The space $(C^\alpha (G) \otimes V_{\sigma_\nu} \otimes \wedge^q(\p_M \oplus \a \oplus \n)^*)^{K_M}$ is the space of $C^\infty$, distribution or hyperfunction sections of the homogeneous vector bundle over $X_M = G/K_M$ defined by the $K_M$-representation on 
$V_{\sigma_\nu}\otimes \wedge^q(\p_M \oplus \a \oplus \n)^*$.
Denote by $\CF^\alpha$ the corresponding sheaves of $C^\infty$-functions, distributions or hyperfunctions so we get a sequence of sheaves
\begin{eqnarray*}
0 \ra (\CF^\alpha (G) \otimes V_{\sigma_\nu})^{K_M} \begin{array}{c}d_1\\ \ra\\ {}\end{array} 
(\CF^\alpha (G) \otimes V_{\sigma_\nu} \otimes \wedge^1(\p_M \oplus \a \oplus \n)^*)^{K_M}\\
\ra (\CF^\alpha (G) \otimes V_{\sigma_\nu} \otimes \wedge^2(\p_M \oplus \a \oplus \n)^*)^{K_M} \ra \dots
\end{eqnarray*}

Let $\K^\alpha$ denote the kernel of $d_1$. We claim that the sequence
\begin{eqnarray*}
0 \ra \K^\alpha \ra (\CF^\alpha (G) \otimes V_{\sigma_\nu})^{K_M} \begin{array}{c}d_1\\ \ra\\ {}\end{array} 
(\CF^\alpha (G) \otimes V_{\sigma_\nu} \otimes \wedge^1(\p_M \oplus \a \oplus \n)^*)^{K_M}\\
\ra (\CF^\alpha (G) \otimes V_{\sigma_\nu} \otimes \wedge^2(\p_M \oplus \a \oplus \n)^*)^{K_M} \ra \dots
\end{eqnarray*}
is exact. By the $G$ equivariance it suffices to prove this at the point $eK_M$ in $X_M$. The exponential map 
$\exp\ : \p_M \oplus \a \oplus \n \oplus \bar{\n} \ra G/K_M$ restricted to a small neighborhood of zero gives local coordinates around the point $eK_M$. The definition of the differential then shows that in those coordinates the differential operators have constant coefficients and are of the type considered in \cite{Kom}, sec. 3.
Thus Theorem 3.2 in \cite{Kom} gives the exactness of the above sequence of sheaves.

The hyperfunction sheaves $\CF^{-\omega}(.)$ are flabby; the sheaves $\CF^\infty$ and $\CF^{-\infty}$ are fine and so this sequence is an acyclic resolution of $\K^\alpha$. 
For the proof of the lemma it remains to show that $\K^\alpha$ is acyclic.
Writing $G=K\times_{K_M} P$ we get 
$X_M =G/K_M = K\times_{K_M}P/M$. This decomposition induces a fibration
$$
P/K_M \ra X_M \= K\times_{K_M}P/K_M \begin{array}{c}f\\ \ra\\ {}\end{array} K/K_M,
$$
with contractible fibre $P/K_M = A\times N \times M/K_M$. The differential $d_1$ is the de Rham differential on the fibres with coefficients in the homogeneous bundle defined by $\sigma |_{K_M}$. So the local sections of $\K^\alpha$ are just those sections of $\CF^\alpha(G\times_{K_M}V_{\sigma_\nu})$ which are locally constant in the fibre direction.
But then we conclude $\K^\alpha = f^{-1}(\CF^\alpha(K\times_{K_M}V_{\sigma_\nu}))$.

Since $f$ has fibre $\cong \R^n$ for some $n$, the extended Vietoris-Begle theorem
\cite{KaSh} Prop. 2.7.8 shows that $H^j(\K^\alpha)=H^j(\CF^\alpha(K\times_{K_M}V_{\sigma_\nu}))$. The latter sheaf is flabby or fine, in any case acyclic, whence the claim.
 \qed

\begin{lemma}
For $\alpha = \infty ,-\infty ,-\omega$ and for any finite dimensional $K_M$-representation $(\tau ,V_\tau)$ the $\Ga$-module $(C^\alpha (G) \otimes V_\tau)^{K_M}\otimes V_\ph$ is $\Ga$-acyclic.
\end{lemma}

\prf
The case $\alpha=-\omega$ is treated analogously to the proof of Lemma 2.6 in \cite{BuOl}. The proof of the other cases follows the lines of the proof of Lemma 2.4 in \cite{BuOl}. \qed

\begin{lemma} \label{gammakohokomplex}
For $\alpha = \infty ,-\infty ,-\omega$ the cohomology of the complex
$$
0 \ra (C^{\alpha}(\Ga \bs G,\ph) \otimes V_{\sigma_\nu})^{K_M}
\ra 
(C^{\alpha}(\Ga \bs G,\ph) \otimes V_{\sigma_\nu} \otimes \wedge^1(\p_M \oplus \a \oplus \n)^*)^{K_M}
$$ $$
\ra (C^{\alpha}(\Ga \bs G,\ph) \otimes V_{\sigma_\nu} \otimes \wedge^2(\p_M \oplus \a \oplus \n)^*)^{K_M} \ra \dots
$$
is isomorphic to $H^*(\Ga ,H_{\sigma ,\nu}^{\alpha} \otimes V_\ph)\cong H^* (\m \oplus \a \oplus \n ,K_M,C^\alpha (\Ga \bs G,\ph))$.
\end{lemma}

\prf
By the last two lemmas the complex (\ref{complex})$\otimes V_\ph$ is a $\Ga$-acyclic resolution of its zeroth cohomology which is just $H_{\sigma ,\nu}^{\alpha}\otimes V_\ph$. \qed

So far we know that the vanishing order at $s=\nu (H_1)$ of $Z_{\theta(P),\sigma ,\ph}(s)$ is
$$
(-1)^{\dim\ \n} \sum_{\pi \in \hat{G}}N_{\Ga ,\ph}(\pi) \sum_{p,q} (-1)^{p+q}
\dim\left(\begin{array}{c}H^q({\n},\pi_K)\otimes \wedge^p\p_M \otimes V_{\sigma}\end{array}\right)^{K_M}_{-\nu}.
$$

For a finite dimensional $K_M$-module $\tau$ and $V$ an arbitrary representation of the group $K_M$ we write $V(\tau)$ for the space $(V\otimes \tau)^{K_M}$.

\begin{lemma} \label{comparison}
For any $\pi \in \hat{G}$ and any finite dimensional representation $\tau$ of $K_M$ we have an isomorphism of finite dimensional $A$-modules
$$
H^q({\n},\pi_K)(\tau) \cong H^q({\n},\pi^\omega)(\tau).
$$
\end{lemma}

\prf
In \cite{Brat} or in \cite{BuOl_conse} it is shown that there is an isomorphism of $A$-modules \\
$H^q(\n ,\pi_K)^{\omega} \cong H^q(\n ,\pi^\omega).$
Therefore the $K$-types agree.\qed

\begin{lemma} \label{duality}
In the situation of the last lemma we have an isomorphism of finite dimensional $A$-modules
$$
H^p({\n},\pi^{-\omega})(\breve{\tau})^* \cong (H^{\dim\ \n -p}({\n},\breve{\pi}^\omega)\otimes \wedge^{\dim\ \n}{\n})(\tau).
$$
\end{lemma}

\prf 
\cite{BuOl} Prop. 4.4.
\qed

Using these lemmas we get the order of $Z_{\theta(P),\sigma ,\ph}(s)$ at $s=\nu(H_1)$ as
$$
\sum_{\pi \in \hat{G}} N_{\Ga ,\breve{\ph}}(\pi) \sum_{p,q}(-1)^{p+q} \dim \left(\begin{array}{c}H^q({\n},\pi^{-\omega})\otimes \wedge^p\p_M\otimes V_{\breve{\sigma}}\end{array}\right)_{-\nu}^{K_M}
$$ $$
= \sum_{p,q} (-1)^{p+q} \dim\left(\begin{array}{c}H^q({\n},C^{-\omega}(\Ga \bs G,\breve{\ph}))\otimes \wedge^p \p_M \otimes V_{\breve{\sigma}}\end{array}\right)_{-\nu}^{K_M}.
$$

\begin{lemma} \label{charpolH}
If $\nu \neq \rho_0$ the $\a$ acts semisimply on $H^q(\m \oplus \n ,K_M,\pi^{-\omega}\otimes V_{\breve{\sigma}})_{-\nu}$. If $\nu =\rho_0$ and $H\in \a$ then $H^2$ still acts semisimply on this space.
\end{lemma}

\prf
The proof is a direct generalization of the corresponding proof in the rank one case (\cite{BuOl} Prop. 4.1). Not to bore the reader with a repetition of the
argumentation in loc.cit. we only indicate the changes: One has to show $(\m \oplus \n ,K_M)$-acyclicity where Bunke and Olbrich show $\n$-acyclicity.
Wherever they use the Iwasawa decomposition one uses the decomposition
$G=NAM\times_{K_M}K$ instead.\qed
In order to prove the theorem let's deal with the case $\nu \neq 0$ first. 
Then, by Lemma \ref{charpolH} the Lie algebra $\a$ acts semisimply on $H^*({\n},C^{-\omega}(\Ga \bs G ,\ph))_{\nu +\rho_0}$ and so we get the order
of $Z_{\theta(P),\sigma ,\ph}(s+|\rho_0|)$ at $s=\nu(H_1)$ as
$$
-\sum_p p(-1)^p \dim\ H^p(\m \oplus \a \oplus \n ,K_M, C^{-\omega}(\Ga \bs G,\breve{\ph})\otimes V_{\breve{\sigma}_{\nu}}),
$$
where Lemma \ref{comparison} and Lemma \ref{duality} assure us that these
cohomology groups are finite dimensional.
By Lemma \ref{gammakohokomplex} this equals
$$
\chi_{_1}(\Ga ,H_{\breve{\sigma},\nu}^{-\omega} \otimes V_{\breve{\ph}}) \= -\sum_p p(-1)^p \dim\ H^p(\Ga ,H_{\breve{\sigma} ,\nu}^{-\omega}\otimes V_{\breve{\ph}}),
$$
where the finite dimensionality follows a fortiori.

Consider the definition of $Z_{P,\sigma ,\ph}$. 
Extending the product over $\ga^{-1}$ instead over $\ga$ we see that $Z_{P,\sigma ,\ph}=Z_{\theta(P),\breve{\sigma},\breve{\ph}}$.
Above we have shown that the order of $Z_{\theta(P),\breve{\sigma},\breve{\ph}}(s)$ at $s=\nu(H_1)$ is
$$
\chi_{_1}(\Ga ,H_{{\sigma},-\nu}^{-\omega} \otimes V_{{\ph}}),
$$
whence the claim.

To finish the proof of the theorem it remains to consider the case $\nu = 0$.
In this case $H\in \a$ does not act semisimply, but $H^2$ does. 
So we consider the complex (\ref{complex}) and change the differential in that we let act $H\in \a$ by $H^2$. The complex remains exact and its zeroth cohomology is a nontrivial extension of $H_{\sigma ,0}^{-\omega}$ by itself.
\qed

\section{The function L} \label{funcL}
In section \ref{constzeta} it was necessary to restrict to locally symmetric manifolds of fundamental rank $\leq 1$. In this section we consider arbitrary compact manifolds $X_\Ga = \Ga \bs X = \Ga \bs G / K$.

Fix a $\theta$-stable Cartan subgroup $H=AB$ of splitrank $r=\dim A$ and a parabolic $P=MAN$. 
Let $\CE_P(\Ga)$ denote the set of all $\Ga$-conjugacy classes $[\ga]$ such that $\ga$ is in $G$ conjugate to an element $a_\ga b_\ga \in A^+B \subset H$.

Let $V$ denote a finite dimensional complex vector space and let $\B = (b_1^*, \dots b_n^*)$ be a basis of the dual space $V^*$ of $V$. The {\bf cross divisor} $D_v$ attached to a point $v\in V$ with respect to the basis $\B$ is by definition
$$
D_v \= \sum_{j=1}^n \{ b_j^* \= b_j^*(v)\} .
$$
A meromorphic function $f$ on $V$ is said to have a simple pole along $D_v$ of residue $s\in \C$ if the function $z \mapsto f(z) - \frac{c}{(b_1^*(z)-b_1^*(v))\dots (b_n^*(z)-b_n^*(v))}$ is holomorphic in an neighborhood of $D_v$.

Fix a basis $\CA = (\alpha_1 ,\dots , \alpha_r)$ of $\a^*$ such that each $\alpha_i$ is a positive multiple of a simple root with respect to $P$ and the volume form $\alpha_1 \dots \alpha_r$ on $\a_0$ coincides with the one induced by the chosen invariant form $B$.

The {\bf flat generated by a closed geodesic} $c$ is the intersection of all maximal compact flats containing $c$. The dimension of this flat is the {\bf rank} of $c$.
Let $\la_c$ denote the volume of this flat. It turns out that on locally symmetric manifolds it holds $\la_c = \la_{c'}$ if $c$ and $c'$ are freely homotopic. So it makes sense to write $\la_{[\ga]}$ or $\la_\ga$.

A {\bf virtual vector space} $V$ is a formal difference $V=V^+-V^-$ of vector spaces. 
We write $\dim^*V=\dim V^+ -\dim V^-$. If $V+\oplus_kV^K$ has a $\Z$-gradation we define $V^+=\oplus_kV^{2k}$ and $V^-=\oplus_kV^{2k+1}$.

Fix a finite dimensional representation $(\tau ,V_\tau)$ of $K_M = K\cap M$.

\begin{theorem}
Let $\Ga$ be neat and $(\ph ,V_\ph)$ a finite dimensional unitary representation of $\Ga$. For $s\in \a^*$ with $\Re (s)>>0$ in the order given by the cone $\a_0^+$ and for $j\in \N$ sufficiently large the Dirichlet series
$$
L^j_{P,\tau ,\ph}(s) \= \sum_{[\ga] \in \CE_P(\Ga)} \la_\ga \chi_{_r}(X_\ga)\ \tr\ \ph(\ga)\ \tr\ \tau(b_\ga)\ \det (1-a_\ga |\n)^j\ \ a^{-s}
$$
converges and extends to a meromorphic function. The poles of $L^j_{P,\tau ,\ph}$ are simple poles along the divisors $D_\la$, $\la \in \a^*$ with respect to the basis $\CA$. The residue at $\la \in \a^*$ is
$$
\sum_{\pi \in \hat{G}} N_{\Ga ,\ph}(^\theta\pi) \dim^* \left( (H^*({\n},\pi_K)\otimes \wedge^*\p_M \otimes \wedge^*\bar{\n} \otimes V_{\tau})^{K_M}\otimes (\wedge^*\bar{\n})^{\otimes j}\right)_{-(\la +2\rho_0)}.
$$
\end{theorem}

\prf
The proof is similar to the proof of Theorem \ref{genSelberg}. 
Let $A^+ \subset A$ denote the positive Weyl chamber given by 
the choice of $P$. Define for $s\in \a^*$, $j\in \N$ the 
function $g_s^{j}$ on $\overline{A^+}$ by $g_s^{j}(a) = 
\det (1-a|\n)^j a^{-s}$.
Choose $\eta : N \rightarrow [0,\infty[$ smooth, $K_M$-invariant, 
with compact support, and such that $\int_N \eta(n) dn =1$. 
Then define
$\tilde{\Phi}(knma(kn)^{-1}) = \eta(n) f_\tau(m)g_s^{j}(a).$

The geometric side of the trace formula is
$$
\tr \left( \tilde{\Phi} | L^2(\Ga \bs G ,\ph)\right) \= 
\sum_{[\ga] \in \CE_P(\Ga)} \la_\ga \chi_{_r}(X_\ga)\ \tr\ \ph(\ga)\ \tr\ \tau(b_\ga)\ \det (1-a_\ga|\n)^j\ a^{-s}.
$$
The spectral side is $\sum_{\pi \in \hat{G}} N_{\Ga ,\ph}(\pi) \tr\ \pi(\tilde{\Phi})$ and we have
$$
\tr\ \pi(\tilde{\Phi}) \= \int_{MA^+} f_\tau(m) \Theta_{H^*(\bar{\n},\pi_K)\otimes \wedge^*\n}(ma) \det (1-a|\n)^j a^{-s-2\rho_0} dm da.
$$
From this the theorem follows.
\qed

\begin{conjecture}The above theorem also holds for small $j$, especially for the case $j=0$.
\end{conjecture}

%\begin{references}

%--------------------Here the manuscript ends-------------------------------
\Addresses
\end{document}